\documentclass[10pt]{iopart}
\usepackage{graphicx,bm,amssymb}
\newcommand{\beq}{\begin{equation}}
\newcommand{\eeq}{\end{equation}}
\newcommand{\beqa}{\begin{eqnarray}}
\newcommand{\eeqa}{\end{eqnarray}}

\newcommand{\abs}[1]{{\left\vert#1\right\vert}}
\newcommand{\ii}{{\rm i}}
\newcommand{\mod}{\mathop{\rm mod}}
\newcommand{\pr}{\mathop{\rm prime}}
\newcommand{\zed}{{\mathbb Z}}
\renewcommand{\Re}{\mathop{\rm Re}}
\newcommand{\Q}{\mathit{\Pi}}
\newcommand{\Z}{{\bm Z}}

\begin{document}

\title{Regularized products of Gauss and Eisenstein integers and primes}

\author{P L Krapivsky$^{1,2}$ and J M Luck$^3$}%\footnote{Corresponding author}

\address{$^1$ Department of Physics, Boston University, Boston, MA 02215, USA}

\address{$^2$ Santa Fe Institute, Santa Fe, NM 87501, USA}

\address{$^3$ Universit\'e Paris-Saclay, CNRS, CEA, Institut de Physique Th\'eorique,
91191~Gif-sur-Yvette, France}

\begin{abstract}
We provide heuristic computations \`a la Euler of the regularized infinite products
of Gauss and Eisenstein integers and primes.
Our approach, yielding explicit expressions,
is inspired by the work by Muñoz Garc\'ia and P\'erez-Marco,
who evaluated the product of all natural primes to $4\pi^2$.
\end{abstract}

\eads{\mailto{pkrapivsky@gmail.com},\mailto{jean-marc.luck@ipht.fr}}

\maketitle

\section{Introduction and summary of results}

Infinite products of integers are obviously divergent.
A broad class of such products may however be given a meaning,
i.e., attributed a well-defined finite value, by means of the zeta-regularization approach.
This technique dates back to Euler, long before it was given a firm basis in terms of analytical continuation
(see~\cite{davis,edwards,vara,varabook} for historical accounts,
and e.g.~\cite{RS,H,V,SABK,QHS,manin,I,AZ} for modern references).
In the physics literature,
zeta-regularized infinite products are often referred to as renormalized products.
A prominent instance consists in the determinants, i.e., products of eigenvalues,
of various operators arising in quantum mechanics and quantum field theory
(see~\cite{EOR} for an overview).

The present note is directly inspired by the work by Muñoz Garc\'ia and P\'erez-Marco~\cite{MP1,MP2},
summarized in section~\ref{natural},
where the product of natural primes was evaluated to
\beq
\Q=4\pi^2.
\eeq
This outcome may be put in perspective with the famous Euler formula for the product of all natural integers,
\beq
P=\sqrt{2\pi}.
\eeq
We thus have
\beq
\Q=P^4.
\label{qpnat}
\eeq

Hereafter we provide heuristic computations \`a la Euler of the regularized infinite products
of Gauss and Eisenstein integers, forming the rings $\zed[\ii]$ and $\zed[\omega]$,
with $\omega=\e^{2\pi\ii/3}$, and of the associated primes.
These examples are the most evident generalizations of the natural integers and primes.
They are also related to the simplest quantum billiards whose spectra are integrable,
namely the square and the equilateral triangle~\cite{IL,I2}.
Our main outcomes are as follows.
For Gauss integers (see section~\ref{gauss}),
we recover that the product $P_4$ of all non-zero Gauss integers is such that
\beq
\abs{P_4}=\frac{\Gamma(1/4)^2}{2\sqrt\pi}\approx3.708149,
\label{p4}
\eeq
and find that the product $\Q_4$ over Gauss primes obeys
\beq
\abs{\Q_4}=\abs{P_4}^8.
\label{qp4}
\eeq
For Eisenstein integers (see section~\ref{eisen}),
we find that the product $P_3$ of all non-zero Eisenstein integers is such that
\beq
\abs{P_3}=\frac{3^{1/4}\,\Gamma(1/3)^3}{2\pi}\approx4.027065,
\label{p3}
\eeq
and that the product $\Q_3$ over Eisenstein primes obeys
\beq
\abs{\Q_3}=\abs{P_3}^{12}.
\label{qp3}
\eeq
The striking resemblance between~(\ref{qpnat}),~(\ref{qp4}) and~(\ref{qp3})
suggests the existence of a more general underlying property, of which we are unaware.

\newpage

\section{A reminder on the products of natural integers and primes}
\label{natural}

We begin with a reminder on the products of natural (positive) integers and primes.

\subsection{Product of natural integers}

The product of all natural integers,
\beq
P=\prod_{n\ge1}n,
\label{pnatural}
\eeq
can be regularized by relating it to the Riemann zeta function,
\beq
\zeta(s)=\sum_{n\ge1}n^{-s}=\prod_{p\pr}(1-p^{-s})^{-1}.
\label{zeta}
\eeq
The identity between both rightmost sides is known as the Euler product formula.
It expresses that every natural integer $n$ can be decomposed in a unique way as a product of primes.
We have formally
\beq
P=\exp\sum_{n\ge1}\ln n=\e^{-\zeta'(0)}.
\label{lnpnatural}
\eeq
This expression is the gist of the zeta-regularization,
consisting in regularizing divergent products such as~(\ref{pnatural})
by means of analytic continuation.
The series and product entering~(\ref{zeta}) converge for $\Re s>1$
and can be analytically continued to a meromorphic function in the whole complex $s$-plane,
with a simple pole with unit residue at $s=1$.
In particular, $\zeta(s)$ is analytic in a neighborhood of $s=0$.
We have
\beq
\zeta(0)=-\frac12,\qquad\zeta'(0)=-\frac12\ln 2\pi.
\label{z0}
\eeq
Inserting the second of these expressions into~(\ref{lnpnatural}), we obtain the famous Euler formula
\beq
P=\prod_{n\ge1}n=\infty\,!=\sqrt{2\pi},
\label{euler}
\eeq
sometimes also attributed to Riemann.

\subsection{Product of natural primes}
\label{naturalprimes}

The product of all natural primes,
\beq
\Q=\prod_{p\pr}p,
\label{qnatural}
\eeq
has been evaluated by Muñoz Garc\'ia and P\'erez-Marco,
first by a heuristic computation \`a la Euler~\cite{MP1}, and then rigorously~\cite{MP2}.
These investigations were motivated by a question raised in~\cite[p.~101]{SABK}.
Remarkably enough, their approach circumvents the obstruction that the prime zeta function $Z(s)$,
to be introduced in~(\ref{zdef}),
has a natural boundary along the whole imaginary axis~\cite{LW},
and is therefore not analytic in a neighborhood of $s=0$.

The heuristic derivation exposed in~\cite{MP1}
starts with the Artin-Hasse identity for the exponential function:
\beq
\e^x=\prod_{n\ge1}(1-x^n)^{-\mu(n)/n},
\label{ah}
\eeq
where $\mu$ is the M\"obius function:
\beq
\mu(n)=\left\{
\begin{array}{cl}
(-1)^k & \hbox{if $n$ is the product of $k$ distinct primes},\\
0 &\hbox{else}.
\end{array}
\right.%}
\eeq
The prime zeta function,
\beq
Z(s)=\sum_{p\pr}p^{-s},
\label{zdef}
\eeq
such that
\beq
\Q=\e^{-Z'(0)},
\label{qz}
\eeq
therefore obeys
\beqa
\e^{Z(s)}
&=&\prod_{p\pr}e^{p^{-s}}
\nonumber\\
&=&\prod_{p\pr}\prod_{n\ge1}(1-p^{-ns})^{-\mu(n)/n}
\nonumber\\
&=&\prod_{n\ge1}\prod_{p\pr}(1-p^{-ns})^{-\mu(n)/n}
\nonumber\\
&=&\prod_{n\ge1}\zeta(ns)^{\mu(n)/n}.
\label{ezprod}
\eeqa
It is underlined in~\cite{MP1,MP2} that the above expression has been known for long~\cite{LW,D}.
Taking logarithmic derivatives of both sides in~(\ref{ezprod}) yields
\beq
Z'(s)=\sum_{n\ge1}\mu(n)\frac{\zeta'(ns)}{\zeta(ns)},
\eeq
and in particular
\beq
Z'(0)=\frac{\zeta'(0)}{\zeta(0)}\sum_{n\ge1}\mu(n).
\label{zpz}
\eeq
The sum entering the above expression is divergent,
but it can consistently be evaluated by means of zeta-regularization, using the identity
\beq
\sum_{n\ge1}\mu(n)\,n^{-s}=\frac{1}{\zeta(s)}.
\eeq
We thus obtain, using~(\ref{z0}),
\beq
Z'(0)=\frac{\zeta'(0)}{\zeta(0)^2}=-2\ln(2\pi),
\label{zp}
\eeq
and so, using~(\ref{qz}),
\beq
\Q=4\pi^2.
\eeq
This is the central result of~\cite{MP1,MP2}.
Notice the identity
\beq
\Q=P^4,
\eeq
announced in~(\ref{qpnat}).

\section{The products of Gauss integers and primes}
\label{gauss}

Gauss integers form the ring $\zed[\ii]$.
They can be viewed as the vertices of the unit square lattice,
of the form $z=a+b\ii$, where $a$ and $b$ are usual integers, and so $\abs{z}^2=a^2+b^2$.

\subsection{Product of Gauss integers}

The product of all Gauss integers is
\beq
P_4=\prod_{(a,b)\ne(0,0)}(a+b\ii).
\label{p4def}
\eeq
The ring $\zed[\ii]$ is invariant under multiplication by $\ii$,
every non-zero Gauss integer having an orbit of length 4:
$(a,b)\to(-b,a)\to(-a,-b)\to(b,-a)$.
In particular, $\zed[\ii]$ has 4 units,
the distinct powers of $\ii$,
namely $1$, $\ii$, $-1$ and $-\ii$.
It therefore seems hopeless to give a meaning to the phase of $P_4$.
Henceforth we rather consider
\beq
\abs{P_4}^2=\prod_{(a,b)\ne(0,0)}(a^2+b^2).
\label{p42def}
\eeq
The corresponding Dedekind zeta function,
\beq
\zeta_4(s)=\sum_{(a,b)\ne(0,0)}(a^2+b^2)^{-s},
\label{z42def}
\eeq
is such that
\beq
\abs{P_4}^2=\e^{-\zeta_4'(0)}.
\label{p42}
\eeq
The theory of Dedekind zeta functions is exposed in detail in~\cite{cartier,kks2,kks3}
for the cases of present interest, namely $\zed[\ii]$ and $\zed[\omega]$,
and in~\cite[Sec.~10]{cohen} for the rings of integers of more general number fields
(see also~\cite{HW} for an overview).
In the present situation, we~have
\beq
\zeta_4(s)=4\zeta(s)L_4(s),
\label{z4factor}
\eeq
where $\zeta(s)$ is the Riemann zeta function,
and
\beq
L_4(s)=\sum_{n\ge1}\chi_4(n)\,n^{-s}
\eeq
is the Dirichlet $L$-function associated with the Dirichlet character $\chi_4$,
i..e, the completely multiplicative function of $n \mod 4$ such that
\beq
\chi_4(0)=0,\quad
\chi_4(1)=1,\quad
\chi_4(2)=0,\quad
\chi_4(3)=-1.
\label{chi4def}
\eeq
This $L$-function reads explicitly
\beqa
L_4(s)
&=&\sum_{m\ge0}\left((4m+1)^{-s}-(4m+3)^{-s}\right)
\nonumber\\
&=&\frac{\zeta(s,1/4)-\zeta(s,3/4)}{4^s},
\label{l4h}
\eeqa
in terms of the Hurwitz zeta function
\beq
\zeta(s,x)=\sum_{m\ge0}(m+x)^{-s}.
\label{hurwitz}
\eeq
The analytic structure of the latter function is similar to that of the Riemann zeta function,
which is recovered as $\zeta(s)=\zeta(s,1)$.
We have in particular
\beq
\zeta(0,x)=\frac12-x,\qquad\zeta'(0,x)=\ln\frac{\Gamma(x)}{\sqrt{2\pi}},
\label{z0x}
\eeq
the second expression being known as the Lerch formula~\cite{lerch}.
We also mention that $L_4(s)$ has an Euler product representation of the form
\beq
L_4(s)=\prod_{p\pr\ne2}(1-\chi_4(p)\,p^{-s})^{-1},
\eeq
i.e., explicitly
\beq
L_4(s)=\frac{\zeta_{4,1}(s)\zeta_{4,3}(2s)}{\zeta_{4,3}(s)},
\label{l4zz}
\eeq
in terms of the partial prime zeta functions
\beqa
\zeta_{4,1}(s)&=&\prod_{p\pr=1\mod4}(1-p^{-s})^{-1},
\nonumber\\
\zeta_{4,3}(s)&=&\prod_{p\pr=3\mod4}(1-p^{-s})^{-1},
\label{z4143def}
\eeqa
which obey
\beq
\zeta_{4,1}(s)\zeta_{4,3}(s)=(1-2^{-s})\zeta(s).
\label{z4143z}
\eeq
Using~(\ref{z4factor}) and~(\ref{l4h}), together with~(\ref{z0}) and~(\ref{z0x}), we obtain
\beqa
&&L_4(0)=\frac12,\qquad
L'_4(0)=\ln\frac{\Gamma(1/4)}{2\Gamma(3/4)},
\label{l4res}
\\
&&\zeta_4(0)=-1,\qquad
\zeta'_4(0)=\ln\frac{2\Gamma(3/4)^2}{\pi\Gamma(1/4)^2}.
\label{zeta4res}
\eeqa
The value of $\zeta_4(0)$ can be interpreted as counting negatively
the point $(0,0)$ that is excluded from the product~(\ref{p4def}) and the sum~(\ref{z42def}).
Finally,~(\ref{p42}) yields
\beq
\abs{P_4}^2=\frac{\pi\Gamma(1/4)^2}{2\Gamma(3/4)^2}=\frac{\Gamma(1/4)^4}{4\pi},
\eeq
and so
\beq
\abs{P_4}=\frac{\Gamma(1/4)^2}{2\sqrt\pi}\approx3.708149.
\eeq
This expression, announced in~(\ref{p4}), can be found e.g.~in~\cite{kks3}.

\subsection{Product of Gauss primes}

The product of all (unnormalized) Gauss primes,
\beq
\Q_4=\prod_{a+b\ii\pr}(a+b\ii),
\label{q4}
\eeq
is also invariant under multiplication by $\ii$,
so that we rather consider
\beq
\abs{\Q_4}^2=\prod_{a+b\ii\pr}(a^2+b^2).
\label{q42def}
\eeq
This quantity can be evaluated along the lines of~\cite{MP1,MP2},
as recalled above,
up to the replacement of the prime zeta function $Z(s)$ by
\beq
Z_4(s)=\sum_{a+b\ii\pr}(a^2+b^2)^{-s},
\eeq
such that
\beq
\abs{\Q_4}^2=\e^{-Z_4'(0)}.
\eeq
Accordingly, in~(\ref{ezprod}) to~(\ref{zpz}), the Riemann zeta function is to be replaced by
\beq
\Z_4(s)=\sum_{a+b\ii\pr}(1-(a^2+b^2)^{-s})^{-1}.
\eeq
This function can be derived from the explicit knowledge of Gauss primes, dating back to Gauss himself.
The Gauss integer $z=a+b\ii$ is prime in the following three situations (see e.g.~\cite[Ch.~XV]{HW}):

\begin{itemize}

\item $z$ is the product of a prime $p=3 \mod 4$ by a unit.
The natural prime $p$ is said to be inert in $\zed[\ii]$.
It corresponds to 4 distinct Gauss primes, as there are 4 units.

\item $\abs{z}^2=a^2+b^2$ is a prime $p=1 \mod 4$.
The natural prime $p$ is said to be split.
It corresponds to 8 distinct Gauss primes,
namely the products of $a+\ii b$ and $b+\ii a$ by units.

\item $z$ is the product of $1+\ii$ by a unit.
The natural prime $\abs{z}^2=2$ is said to be ramified.
It corresponds to 4 distinct Gauss primes.

\end{itemize}

Taking the above multiplicities into account, and using~(\ref{z4143def}), we obtain
\beq
\Z_4(s)=\zeta_{4,3}(2s)^4\zeta_{4,1}(s)^8(1-2^{-s})^{-4}.
\eeq
Using~(\ref{l4zz}) and~(\ref{z4143z}), this boils down to
\beq
\Z_4(s)=\zeta(s)^4L_4(s)^4=\left(\frac{\zeta_4(s)}{4}\right)^4.
\eeq
The analogue of~(\ref{zp}) therefore reads
\beq
Z_4'(0)
=\frac{1}{\zeta(0)}\,\frac{\Z_4'(0)}{\Z_4(0)}
=\frac{4\zeta_4'(0)}{\zeta(0)\zeta_4(0)}.
\eeq
Using~(\ref{z0}) and~(\ref{zeta4res}), we obtain
\beq
Z_4'(0)=8\zeta_4'(0).
\eeq
We are thus left with the result
\beq
\abs{\Q_4}=\abs{P_4}^8,
\eeq
announced in~(\ref{qp4}).

\newpage

\section{The products of Eisenstein integers and primes}
\label{eisen}

Eisenstein integers form the ring $\zed[\omega]$.
They can be viewed as the vertices of the unit triangular lattice,
of the form $z=a+b\omega$, where $a$ and $b$ are usual integers,
and $\omega=\e^{2\pi\ii/3}$, and so $\abs{z}^2=a^2-ab+b^2$.

\subsection{Product of Eisenstein integers}

The product of all Eisenstein integers is
\beq
P_3=\prod_{(a,b)\ne(0,0)}(a+b\omega).
\label{p3def}
\eeq
The ring $\zed[\omega]$ is invariant under multiplication by $\rho=\e^{\ii\pi/3}=\sqrt{\omega}=1+\omega$,
every non-zero Eisenstein integer having an orbit of length 6:
$(a,b)\to(a-b,a)\to(-b,a-b)\to(-a,-b)\to(b-a,-a)\to(b,b-a)$.
In particular, $\zed[\omega]$ has 6 units,
the distinct powers of $\rho$,
namely $1$, $1+\omega$, $\omega$, $-1$, $-1-\omega$ and $-\omega$.
It again seems hopeless to give a meaning to the phase of $P_3$.
We rather consider
\beq
\abs{P_3}^2=\prod_{(a,b)\ne(0,0)}(a^2-ab+b^2).
\label{p32def}
\eeq
The corresponding Dedekind zeta function,
\beq
\zeta_3(s)=\sum_{(a,b)\ne(0,0)}(a^2-ab+b^2)^{-s},
\label{z32def}
\eeq
such that
\beq
\abs{P_3}^2=\e^{-\zeta_3'(0)},
\label{p32}
\eeq
reads
\beq
\zeta_3(s)=6\zeta(s)L_3(s),
\label{z3factor}
\eeq
where
\beq
L_3(s)=\sum_{n\ge1}\chi_3(n)\,n^{-s}
\eeq
is the Dirichlet $L$-function associated with the Dirichlet character $\chi_3$,
i..e, the completely multiplicative function of $n \mod 3$ such that
\beq
\chi_3(0)=0,\quad
\chi_3(1)=1,\quad
\chi_3(2)=-1.
\label{chi3def}
\eeq
This $L$-function reads explicitly
\beqa
L_3(s)
&=&\sum_{m\ge0}\left((3m+1)^{-s}-(3m+2)^{-s}\right)
\nonumber\\
&=&\frac{\zeta(s,1/3)-\zeta(s,2/3)}{3^s}.
\label{l3h}
\eeqa
We also mention that $L_3(s)$ has an Euler product representation of the form
\beq
L_3(s)=\prod_{p\pr\ne3}(1-\chi_3(p)\,p^{-s})^{-1},
\eeq
i.e., explicitly
\beq
L_3(s)=\frac{\zeta_{3,1}(s)\zeta_{3,2}(2s)}{\zeta_{3,2}(s)},
\label{l3zz}
\eeq
in terms of the partial prime zeta functions
\beqa
\zeta_{3,1}(s)&=&\prod_{p\pr=1\mod3}(1-p^{-s})^{-1},
\nonumber\\
\zeta_{3,2}(s)&=&\prod_{p\pr=2\mod3}(1-p^{-s})^{-1},
\label{z3132def}
\eeqa
which obey
\beq
\zeta_{3,1}(s)\zeta_{3,2}(s)=(1-3^{-s})\zeta(s).
\label{z3132z}
\eeq
Using~(\ref{z3factor}) and~(\ref{l3h}), together with~(\ref{z0}) and~(\ref{z0x}), we obtain
\beqa
&&L_3(0)=\frac13,\qquad
L'_3(0)=\ln\frac{\Gamma(1/3)}{3^{1/3}\Gamma(2/3)},
\label{l3res}
\\
&&\zeta_3(0)=-1,\qquad
\zeta'_3(0)=\ln\frac{3\Gamma(2/3)^3}{2\pi\Gamma(1/3)^3}.
\label{zeta3res}
\eeqa
The value of $\zeta_3(0)$ can again be interpreted as counting negatively
the point $(0,0)$ that is excluded from the product~(\ref{p3def}) and the sum~(\ref{z32def}).
Finally,~(\ref{p32}) yields
\beq
\abs{P_3}^2=\frac{2\pi\Gamma(1/3)^3}{3\Gamma(2/3)^3}=\frac{\sqrt{3}\,\Gamma(1/3)^6}{4\pi^2},
\eeq
and so
\beq
\abs{P_3}=\frac{3^{1/4}\,\Gamma(1/3)^3}{2\pi}\approx4.027065.
\eeq
This expression is announced in~(\ref{p3}).

\subsection{Product of Eisenstein primes}

The product of all (unnormalized) Eisenstein primes,
\beq
\Q_3=\prod_{a+b\omega\pr}(a+b\omega),
\label{q3}
\eeq
is also invariant under multiplication by $\rho$,
so that we rather consider
\beq
\abs{\Q_3}^2=\prod_{a+b\omega\pr}(a^2-ab+b^2).
\label{q32def}
\eeq
This quantity can again be evaluated along the lines of~\cite{MP1,MP2},
up to the replacement of the prime zeta function $Z(s)$ by
\beq
Z_3(s)=\sum_{a+b\omega\pr}(a^2-ab+b^2)^{-s},
\eeq
such that
\beq
\abs{\Q_3}^2=\e^{-Z_3'(0)}.
\eeq
Accordingly, in~(\ref{ezprod}) to~(\ref{zpz}), the Riemann zeta function is to be replaced by
\beq
\Z_3(s)=\sum_{a+b\omega\pr}(1-(a^2-ab+b^2)^{-s})^{-1}.
\eeq
This function can be derived from the explicit knowledge of Eisenstein primes.
The Eisenstein integer $z=a+b\omega$ is prime in the following three situations (see e.g.~\cite[Ch.~XV]{HW}):

\begin{itemize}

\item $z$ is the product of a prime $p=2 \mod 3$ by a unit.
The natural prime $p$ is said to be inert in $\zed[\omega]$.
It corresponds to 6 distinct Eisenstein primes, as there are 6 units.

\item $\abs{z}^2=a^2-ab+b^2$ is a prime $p=1 \mod 3$.
The natural prime $p$ is said to be split.
It corresponds to 12 distinct Eisenstein primes.

\item $z$ is the product of $1-\omega$ by a unit.
The natural prime $\abs{z}^2=3$ is said to be ramified.
It corresponds to 6 distinct Eisenstein primes.

\end{itemize}

Taking the above multiplicities into account, and using~(\ref{z3132def}), we obtain
\beq
\Z_3(s)=\zeta_{3,2}(2s)^6\zeta_{3,1}(s)^{12}(1-3^{-s})^{-6}.
\eeq
Using~(\ref{l3zz}) and~(\ref{z3132z}), this boils down to
\beq
\Z_3(s)=\zeta(s)^6L_3(s)^6=\left(\frac{\zeta_3(s)}{6}\right)^6.
\eeq
The analogue of~(\ref{zp}) therefore reads
\beq
Z_3'(0)
=\frac{1}{\zeta(0)}\,\frac{\Z_3'(0)}{\Z_3(0)}
=\frac{6\zeta_3'(0)}{\zeta(0)\zeta_3(0)}.
\eeq
Using~(\ref{z0}) and~(\ref{zeta3res}), we obtain
\beq
Z_3'(0)=12\zeta_3'(0).
\eeq
We are thus left with the result
\beq
\abs{\Q_3}=\abs{P_3}^{12},
\eeq
announced in~(\ref{qp3}).

\newpage

\section*{References}

\bibliography{products.bib}

\end{document}